\documentclass[12pt,twoside]{article}
\usepackage[english]{babel}
\usepackage[latin1]{inputenc}
\usepackage{amsmath}
\usepackage{mathtools} 
\usepackage{graphicx}
\usepackage{times,amssymb,amscd}
\newtheorem{thm}{Theorem}[section]

\newtheorem{conj}[thm]{Conjecture}
\newtheorem{lem}[thm]{Lemma}

\newtheorem{defn}[thm]{Definition}
\newtheorem{rem}[thm]{Remark}
\numberwithin{equation}{section}

\newcommand{\bA}{\mathbf{A}}

\newcommand{\bE}{\mathbf{E}}

\newcommand{\bH}{\mathbf{H}}

\newcommand{\bL}{\mathbf{L}}

\newcommand{\bR}{\mathbf{R}}
\newcommand{\bS}{\mathbf{S}}
\newcommand{\bV}{\mathbf{V}}

\newcommand{\be}{\mathbf{e}}

\newcommand{\bT}{\mathbf{T}}

\newcommand{\bu}{\mathbf{u}}
\newcommand{\bv}{\mathbf{v}}
\newcommand{\bt}{\mathbf{t}}

\newcommand{\EUC}{\mathbf E^3}
\newcommand{\SPH}{\bS^3}
\newcommand{\HYP}{\bH^3}
\newcommand{\SXR}{\bS^2\!\times\!\bR}
\newcommand{\HXR}{\bH^2\!\times\!\bR}
\newcommand{\SLR}{\widetilde{\bS\bL_2\bR}}
\newcommand{\NIL}{\mathbf{Nil}}
\newcommand{\SOL}{\mathbf{Sol}}

\begin{document}
\pagestyle{myheadings}
\markboth{\centerline{Jen\H o Szirmai}}
{Triangle angle sums of geodesic triangles and $\dots$}
\title
{Interior angle sums of geodesic triangles and translation-like isoptic surfaces in $\SOL$ geometry
}

\author{G\'eza Csima and Jen\H o Szirmai \\
\normalsize Department of Algebra and Geometry, Institute of Mathematics,\\
\normalsize Budapest University of Technology and Economics, \\
\normalsize M\"uegyetem rkp. 3., H-1111 Budapest, Hungary \\
\normalsize csimageza@gmail.com,~szirmai@math.bme.hu
\date{\normalsize{\today}}}

\maketitle
\begin{abstract}
After having investigated the geodesic triangles and their angle sums in $\NIL$ and $\SLR$ geometries we consider the analogous problem in $\SOL$ space that
is one of the eight 3-dimensional Thurston geometries. We analyse the interior angle sums of geodesic triangles and we prove  
that it can be larger than, less than or equal to $\pi$. 

Moreover, we determine the equations of $\SOL$ isoptic
surfaces of translation-like segments and 
as a special case of this we examine the $\SOL$ translati\-on-like Thales sphere, which we call {\it Thaloid}. We also discuss the behavior of this surface.

In our work we will use the projective model of $\SOL$ described by E. Moln\'ar in \cite{M97}.
\end{abstract}

\section{Introduction} \label{section1}
\subsection{Preliminaries to geodesic and translation triangles}
A geodesic triangle in Riemannian geometry and more generally in metric geometry is a 
figure consisting of three different points together with the pairwise connecting geodesic curves.
The points are called vertices, while the geodesic curves are called sides of the triangle.

In the geometries of constant curvature $\EUC$, $\HYP$, $\SPH$ the well-known sums of the interior angles of geodesic
triangles characterize space. It is related to the Gauss-Bonnet theorem which states that the integral of the Gauss curvature
on a compact $2$-dimensional Riemannian manifold $M$ is equal to $2\pi\chi(M)$ where $\chi(M)$ denotes the Euler characteristic of $M$.
This theorem has a generalization to any compact even-dimensional Riemannian manifold (see e.g. \cite{Ch}, \cite{KN}).

In the Thurston spaces can be introduced in a natural way (see \cite{M97}) translation curves. These curves are simpler than geodesics and 
differ from them in $\NIL$, $\SLR$ and $\SOL$ geometries  (\cite{B}, \cite{Cs-Sz23}, \cite{Sz18}, \cite{Sz20}). In $\EUC$, $\SPH$, $\HYP$, $\SXR$ and $\HXR$ geometries the mentioned curves 
coincide with each other. Similar to the geodesic triangle, we can also introduce the 
notion of the translation triangle in each Thurston geometry.
Previously, in the \cite{CsSz16}, \cite{Sz16}, \cite{Sz202} articles, we examined the sum of the angles of the geodesic triangles of the $\SXR$, $\HXR$, $\NIL$ and $\SLR$ geometries. 
Also, in \cite{CsSz16} we analyzed the possible angle sums of translational triangles in $\SLR$ geometry.

In \cite{Sz20}, we examined the sum of the interior angles of $\SOL$ translational triangles and proved, that it is larger or equal than $\pi$.

In this paper we consider the analogous problem for geodesic triangles in $\SOL$ geometry.
\begin{rem}
We note here, that nowadays the $\SOL$ geometry is a widely investigated space concerning
its manifolds, tilings, geodesic and translation ball packings and probability theory
(see e.g. \cite{BT}, \cite{CaMoSpSz}, \cite{KV}, \cite{MSz}, \cite{MSz12}, \cite{MSzV}, 
\cite{Sz18}, \cite{Sz13-2} and the references given there).
\end{rem}

Now, we are interested in {\it geodesic triangles} in $\SOL$ space that is one of the eight Thurston geometries
\cite{S,T} on the base of Heisenberg matrix group. In Section 2 we describe the projective model of $\SOL$
and we shall use its standard Riemannian metric obtained by
pull back transform to the infinitesimal  arc-length-square at the origin. We also recall the isometry group of $\SOL$ and
give an overview about geodesic curves.

In Section 3 we study the $\SOL$ geodesic triangles and prove that the interior angle sum of a geodesic triangle
in $\SOL$ geometry can be larger than, less than or equal to $\pi$.
\subsection{Preliminaries to isoptic curves}

It is well known that in the Euclidean plane the locus of points from a 
segment subtends a given angle $\alpha$ $(0<\alpha<\pi)$ is the union of 
two arcs except for the endpoints with the segment as common chord. If this $\alpha$ is equal to $\frac{\pi}{2}$ then we get the Thales circle. 
Replacing the segment to another general curve, we obtain the Euclidean definition of 
isoptic curve:
\begin{defn}[\cite{yates}]The locus of the intersection of tangents
to a curve meeting at a constant angle $\alpha$ $(0<\alpha <\pi)$ is
the $\alpha$ -- isoptic of the given curve. The isoptic 
curve with right angle called \textit{orthoptic curve}.
\label{defiso}
\end{defn}
\begin{rem}Sometimes we consider the $\alpha$ -- and $\pi-\alpha$ -- isoptics together. Thus, in the case of the section, we get two circles with the segment as a common chord (endpoints of the segment are excluded). Hereafter, we call them $\alpha$ -- isoptic circles.\end{rem}
Although the name "isoptic curve" was suggested by Taylor  in 1884 (\cite{T}), reference to former results can be found in \cite{yates}. In the obscure history of isoptic curves, we can find the names of la Hire (cycloids 1704) and Chasles (conics and epitrochoids 1837) among the contributors of the subject. A very interesting table of isoptic and orthoptic curves is introduced in \cite{yates}, unfortunately without any exact reference of its source. However, recent works are available on the topic, which shows its timeliness. In \cite{harom} and \cite{tizenketto}, the Euclidean isoptic curves of closed strictly convex curves are studied using their support function.
Papers \cite{Kurusa,Wu71-1,Wu71-2} deal with Euclidean curves having a
circle or an ellipse for an isoptic curve. Further curves appearing as isoptic curves are well studied in Euclidean plane geometry $\mathbf{E}^2$, see e.g. \cite{Loria,Wi}.
Isoptic curves of conic sections have been studied in \cite{H} and \cite{Si}. There are results for Bezier curves by Kunkli et al. as well, see \cite{Kunkli}. Many papers focus on the properties of isoptics, e.g. \cite{nyolc,MM3,het}, and the references  therein. There are some generalizations of the isoptics as well \textit{e.g.} equioptic curves in \cite{Odehnal} by Odehnal or secantopics in \cite{tizenegy, tiz} by Skrzypiec.

The notion of isoptic curve can be extended to the other planes of constant curvature (hyperbolic plane $\mathbf{H}^2$ and spherical plane $\mathbf{H}^2$). We studied these questions in \cite{Cs-Sz14} and \cite{Cs-Sz16-1}.

We have extended the problem of isoptic surfaces for spatial segments to the Euclidean and $\NIL$ space in \cite{Cs-Sz23}. 

For further isoptic surfaces in Euclidean geometry, see \cite{Cs-Sz13,Cs-Sz16}, where we extend the definition of isoptic surfaces to other spatial objects.

\section{On Sol geometry}
\label{sec:1}

In this Section we summarize the significant notions and notations of real $\SOL$ geometry (see \cite{M97}, \cite{S}).

$\SOL$ is defined as a 3-dimensional real Lie group with multiplication  
\begin{equation}
     \begin{gathered}
(a,b,c)(x,y,z)=(x + a e^{-z},y + b e^z ,z + c).
     \end{gathered} \tag{2.1}
     \end{equation}
We note that the conjugacy by $(x,y,z)$ leaves invariant the plane $(a,b,c)$ with fixed $c$:
\begin{equation}
     \begin{gathered}
(x,y,z)^{-1}(a,b,c)(x,y,z)=(x(1-e^{-c})+a e^{-z},y(1-e^c)+b e^z ,c).
     \end{gathered} \tag{2.2}
     \end{equation}
Moreover, for $c=0$, the action of $(x,y,z)$ is only by its $z$-component, where $(x,y,z)^{-1}=(-x e^{z}, -y e^{-z} ,-z)$. Thus the $(a,b,0)$ plane is distinguished as a {\it base plane} in
$\SOL$, or by other words, $(x,y,0)$ is normal subgroup of $\SOL$.
$\SOL$ multiplication can also be affinely (projectively) interpreted by "right translations" 
on its points as the following matrix formula shows, according to (2.1):
     \begin{equation}
     \begin{gathered}
     (1;a,b,c) \to (1;a,b,c)
     \begin{pmatrix}
         1&x&y&z \\
         0&e^{-z}&0&0 \\
         0&0&e^z&0 \\
         0&0&0&1 \\
       \end{pmatrix}
       =(1;x + a e^{-z},y + b e^z ,z + c)
       \end{gathered} \tag{2.3}
     \end{equation}
by row-column multiplication.
This defines "translations" $\mathbf{L}(\mathbf{R})= \{(x,y,z): x,~y,~z\in \mathbf{R} \}$ 
on the points of space $\SOL= \{(a,b,c):a,~b,~c \in \mathbf{R}\}$. 
These translations are not commutative, in general. 
Here we can consider $\mathbf{L}$ as projective collineation group with right actions in homogeneous 
coordinates as usual in classical affine-projective geometry.
We will use the Cartesian homogeneous coordinate simplex $E_0(\be_0)$, $E_1^{\infty}(\be_1)$, \ $E_2^{\infty}(\be_2)$, \ 
$E_3^{\infty}(\be_3), \ (\{\be_i\}\subset \bV^4$ \ $\text{with the unit point}$ $E(\be = \be_0 + \be_1 + \be_2 + \be_3 ))$ 
which is distinguished by an origin $E_0$ and by the ideal points of coordinate axes, respectively. 
Thus {$\SOL$} can be visualized in the affine 3-space $\bA^3$
(so in Euclidean space $\bE^3$) as well \cite{M97}.

In this affine-projective context E. Moln\'ar has derived in \cite{M97} the usual infinitesimal arc-length square at any point 
of $\SOL$, by pull back translation, as follows
\begin{equation}
   \begin{gathered}
      (ds)^2:=e^{2z}(dx)^2 +e^{-2z}(dy)^2 +(dz)^2.
       \end{gathered} \tag{2.4}
     \end{equation}
Hence we get infinitesimal Riemann metric invariant under translations, by the symmetric metric tensor field $g$ on $\SOL$ by components as usual.

It will be important for us that the full isometry group Isom$(\SOL)$ has eight components, since the stabilizer of the origin 
is isomorphic to the dihedral group $\mathbf{D_4}$, generated by two involutive (involutory) transformations, preserving (2.4): 
\begin{equation}
   \begin{gathered}
      (1)  \ \ y \leftrightarrow -y; \ \ (2)  \ x \leftrightarrow y; \ \ z \leftrightarrow -z; \ \ \text{i.e. first by $3\times 3$ matrices}:\\      
     (1) \ \begin{pmatrix}
               1&0&0 \\
               0&-1&0 \\
               0&0&1 \\
     \end{pmatrix}; \ \ \ 
     (2) \ \begin{pmatrix}
               0&1&0 \\
               1&0&0 \\
               0&0&-1 \\
     \end{pmatrix}; \\
     \end{gathered} \tag{2.5}
     \end{equation}
     with its product, generating a cyclic group $\mathbf{C_4}$ of order 4
     \begin{equation}
     \begin{gathered}
     \begin{pmatrix}
                    0&1&0 \\
                    -1&0&0 \\
                    0&0&-1 \\
     \end{pmatrix};\ \ 
     \begin{pmatrix}
               -1&0&0 \\
               0&-1&0 \\
               0&0&1 \\
     \end{pmatrix}; \ \ 
     \begin{pmatrix}
               0&-1&0 \\
               1&0&0 \\
               0&0&-1 \\
     \end{pmatrix};\ \ 
     \mathbf{Id}=\begin{pmatrix}
               1&0&0 \\
               0&1&0 \\
               0&0&1 \\
     \end{pmatrix}. 
     \end{gathered} \notag
     \end{equation}
     Or we write by collineations fixing the origin $O(1,0,0,0)$:
\begin{equation}
(1) \ \begin{pmatrix}
         1&0&0&0 \\
         0&1&0&0 \\
         0&0&-1&0 \\
         0&0&0&1 \\
       \end{pmatrix}, \ \
(2) \ \begin{pmatrix}
         1&0&0&0 \\
         0&0&1&0 \\
         0&1&0&0 \\
         0&0&0&-1 \\
       \end{pmatrix} \ \ \text{of form (2.3)}. \tag{2.6}       
\end{equation}
A general isometry of $\SOL$ to the origin $O$ is defined by a product $\gamma_O \tau_X$, first $\gamma_O$ of form (2.6) then $\tau_X$ of (2.3). To
a general point $A(1,a,b,c)$, this will be a product $\tau_A^{-1} \gamma_O \tau_X$, mapping $A$ into $X(1,x,y,z)$. 

Conjugacy of translation $\tau$ by an above isometry $\gamma$, as $\tau^{\gamma}=\gamma^{-1}\tau\gamma$ also denotes it, will also be used by 
(2.3) and (2.6) or also by coordinates with above conventions. We note here that the {\bf Sol}-space is
translation-complete, i.e. every two points of the {\bf Sol}-space can be
connected through one translation arc and every three points form a triangle, when
they are not on the same translation curve. 

We remark only that the role of $x$ and $y$ can be exchanged throughout the paper, but this leads to the mirror interpretation of $\SOL$.
As formula (2.4) fixes the metric of $\SOL$, the change above is not an isometry of a fixed $\SOL$ interpretation. Other conventions are also accepted
and used in the literature.

{\it $\SOL$ is an affine metric space (affine-projective one in the sense of the unified formulation of \cite{M97}). Therefore its linear, affine, unimodular,
etc. transformations are defined as those of the embedding affine space.}
\subsection{Geodesic curves} \label{subsection2}
The geodesic curves of the $\SOL$ geometry are generally defined as having locally minimal arc length between their any two (near enough) points.
The equation systems of the parameterized geodesic curves $g(x(t),y(t),z(t))$ in our model (now by (2.4)) can be determined by the
Levy-Civita theory of Riemann geometry.
We can assume, that the starting point of a geodesic curve is the origin because we can transform a curve into an
arbitrary starting point by translation (2.1);
The system of differential equation of geodesics can be determined using the usual method of Riemannian geometry 
using our metric (fundamental) tensor, which is derived by (2.4) (see \cite{MSz} and \cite{Sz22-3}).
\begin{equation}
\begin{gathered}
\ddot{x}(t)+2\dot{x}(t) \dot{z}(t) =0, \\
\ddot{y}(t)+2\dot{y}(t) \dot{z}(t) =0, \\
\ddot{z}(t)-e^{2z(t)} (\dot{x}(t))^2 + e^{-2z(t)} (\dot{y}(t))^2.
\end{gathered} \tag{2.7}
\end{equation}
where initial conditions are the following
\begin{equation}
\begin{gathered}
        x(0)=y(0)=z(0)=0; \ \ \dot{x}(0)=u=\cos{\theta} \cdot \cos{\alpha}, \ \dot{y}(0)=v=\cos{\theta} \cdot \sin{\alpha}, \\ \dot{z}(0)=w=\sin{\theta}; \ 
        - \pi \leq \alpha \leq \pi, \ \text{and} \ -\frac{\pi}{2} \leq \theta \leq \frac{\pi}{2} \notag
\end{gathered}
\end{equation}
i.e. unit velocity can be assumed.

Basically, we can distinguish $4$ types of geodesic curves $\gamma(t)=(x(t),y(t),z(t))$: 
\begin{enumerate}
\item $u \cdot v = 0$, $u^2+v^2 \ne 0$, $0 \le |{w}| <1$. 
\begin{equation}
\begin{gathered}
x(t)=0, ~ 
y(t)=u\frac{\sinh{t}}{\cosh{t}-w\cdot \sinh{t}}, \\
z(t)=-\log({\cosh{t}-w\cdot \sinh{t}}) 
\end{gathered} \tag{2.8}
\end{equation}
\begin{rem} We note here that the $\SOL$ isometry (2) described in (2.5) maps the $[y,z]$ 
coordinate plane to the $[x,z]$ coordinate plane, so it is sufficient to 
consider the geodesics of the $[y,z]$ coordinate plane.
\end{rem}
\item $w=0$ and $u \cdot v \ne 0$. 
\begin{equation}
\begin{gathered}
x(t)=vt,~
y(t)=ut,~
z(t)= \frac{1}{2}\log{\Big|\frac{u}{v}\Big|}.
\end{gathered} \tag{2.9}
\end{equation}
This means that in the coordinate plane $[x,y]$ 
for parameters $u=v=\frac{1}{\sqrt{2}}$ the geodesic curve can be derived by the 
straight lines $\gamma(t)=\frac{1}{\sqrt{2}}(\pm t, \pm t, 0)$. 
\item $u=v=0$ and $|{w}|=1$.
\begin{equation}
\begin{gathered}
x(t)=0,~
y(t)=0,~
z(t)= \pm t, \ \text{iff} \ w=\pm 1 
\end{gathered} \tag{2.10}
\end{equation}
\item $u \cdot v \ne 0$, $0 < |{w}| < 1 $.
In these cases 
\begin{equation}
\begin{gathered}
x(t) = u \int_0^t e^{-2z(\tau)} d \tau,~
y(t) = v \int_0^t e^{2z(\tau)} d \tau \\
\text{where} ~ z(t) ~ \text{comes from the separable differential equation} \\
\frac{\pm dz}{\sqrt{1-u^2 e^{-2z}-v^2 e^{2z}}}=dt,~ 
\text{whose solution is non-elementary function.}
\end{gathered} \tag{2.11}
\end{equation}
\end{enumerate}
\begin{rem}
The metric on the plane $[x,z]$ induced by $(ds)^2 =e^{2z}(dx)^2 +(dz)^2$ agrees with the $\SOL$ metric.
On mapping this plane to the upper half plane $\{(x_1,x_2) \in \mathbf{R}:~x_2 > 0$ by mapping 
$\gamma:(x,z) \mapsto (x_1,x_2) =(x,e^{-z})$. We obtain, that the metric becomes 
the standard hyperbolic metric (see Fig.~1). 
\begin{figure}[ht]
\centering
\includegraphics[width=8cm]{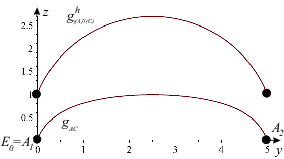}
\caption{The geodesic curve arc $g_{AC}$ lying in the $[y,z]$ plane and its 
image by map $\gamma$, which 
is a circular arc $g^h_{\gamma(A) \gamma(C)}$ in the half-space model of the hyperbolic plane.}
\label{fig:Fig1}
\end{figure}
\begin{equation}
(ds)^2 = \frac{(dx_1)^2+(dx_2)^2}{x_2^2}. \notag
\end{equation}
Therefore, the $[x,z]$ plane (and so the plane $[y,z]$) is a convexly embedded copy 
of $\mathbf{H}^2$ in $\SOL$. Moreover, a similar statement can be formulated for the planes $x=c$, $y=c$ $(c\in \mathbf{R})$.

\end{rem}
\begin{defn}
The distance $d^S(P_1,P_2)$ between the points $P_1$ and $P_2$ lying in the $\SOL$ 
space is defined by the arc length of geodesic curve from $P_1$ to $P_2$.
\end{defn}
\subsection{Translation curves}

We consider a $\SOL$ curve $(1,x(t), y(t), z(t) )$ with a given starting tangent vector at the origin $O(1,0,0,0)$
\begin{equation}
   \begin{gathered}
      u=\dot{x}(0),\ v=\dot{y}(0), \ w=\dot{z}(0).
       \end{gathered} \tag{2.12}
     \end{equation}
For a translation curve let its tangent vector at the point $(1,x(t), y(t), z(t) )$ be defined by the matrix (2.3) 
with the following equation:
\begin{equation}
     \begin{gathered}
     (0,u,v,w)
     \begin{pmatrix}
         1&x(t)&y(t)&z(t) \\
         0&e^{-z(t)}&0&0 \\
         0&0&e^{z(t)}& 0 \\
         0&0&0&1 \\
       \end{pmatrix}
       =(0,\dot{x}(t),\dot{y}(t),\dot{z}(t)).
       \end{gathered} \tag{2.13}
     \end{equation}
Thus, {\it translation curves} in $\SOL$ geometry (see \cite{MoSzi10} and \cite{MSz}) are defined by the first order differential equation system 
$\dot{x}(t)=u e^{-z(t)}, \ \dot{y}(t)=v e^{z(t)},  \ \dot{z}(t)=w,$ whose solution is the following: 
\begin{equation}
   \begin{gathered}
     x(t)=-\frac{u}{w} (e^{-wt}-1), \ y(t)=\frac{v}{w} (e^{wt}-1),  \ z(t)=wt, \ \mathrm{if} \ w \ne 0 \ \mathrm{and} \\
     x(t)=u t, \ y(t)=v t,  \ z(t)=z(0)=0 \ \ \mathrm{if} \ w =0.
       \end{gathered} \tag{2.14}
\end{equation}

We assume that the starting point of a translation curve is the origin, because we can transform a curve into an 
arbitrary starting point by translation (2.3), moreover, unit velocity translation can be assumed :
\begin{equation}
\begin{gathered}
        x(0)=y(0)=z(0)=0; \\ \ u=\dot{x}(0)=\cos{\theta} \cos{\phi}, \ \ v=\dot{y}(0)=\cos{\theta} \sin{\phi}, \ \ w=\dot{z}(0)=\sin{\theta}; \\ 
        - \pi \leq \phi \leq \pi, \ -\frac{\pi}{2} \leq \theta \leq \frac{\pi}{2}. \tag{2.15}
\end{gathered}
\end{equation}
\begin{defn}
The translation distance $d^t(P_1,P_2)$ between the points $P_1$ and $P_2$ is defined by the arc length of the above translation curve 
from $P_1$ to $P_2$.
\end{defn}
Thus we obtain the parametric equation of the the {\it translation curve segment} $t(\phi,\theta,t)$ with starting point at the origin in direction 
\begin{equation}
\bt(\phi, \theta)=(\cos{\theta} \cos{\phi}, \cos{\theta} \sin{\phi}, \sin{\theta}) \tag{2.16}
\end{equation}
where $t \in [0,r \in \bR^+$]. If $\theta \ne 0$ then the system of equation is:
\begin{equation}
\begin{gathered}
        \left\{ \begin{array}{ll} 
        x(\phi,\theta,t)=-\cot{\theta} \cos{\phi} (e^{-t \sin{\theta}}-1), \\
        y(\phi,\theta,t)=\cot{\theta} \sin{\phi} (e^{t \sin{\theta}}-1), \\
        z(\phi,\theta,t)=t \sin{\theta}.
        \end{array} \right. \\
        \text{If $\theta=0$ then}: ~  x(\phi,\theta,t)=t\cos{\phi} , \ y(\phi,\theta,t)=t \sin{\phi},  \ z(\phi,\theta,t)=0.
        \tag{2.17}
\end{gathered}
\end{equation}
It will be important for us later on to be able to invert the mapping in (2.17) to obtain the value of $\phi$ and $\theta$ angles for a given point with homogeneous coordinates $(1,x,y,z)$ (\cite{Sz18}, \cite{Sz20} ):
\begin{lem}
\begin{enumerate}
\item Let $(1,x,y,z)$ $(y,z \in \bR \setminus \{0\}, x \in \bR)$ be the homogeneous coordinates of the point $P \in \SOL$. The parameters of the 
corresponding translation curve $t_{E_0P}$ are the following
\begin{equation}
\begin{gathered}
\phi=\mathrm{arccot}\Big(-\frac{x}{y} \frac{\mathrm{e}^z-1}{\mathrm{e}^{-z}-1}\Big),~\theta=\mathrm{arccot}\Big( \frac{y}{\sin\phi(\mathrm{e}^z-1)}\Big),\\
t=\frac{z}{\sin\theta}, ~ \text{where} ~ -\pi < \phi \le \pi, ~ -\pi/2\le \theta \le \pi/2, ~ t\in \bR^+.
\end{gathered} \tag{2.18}
\end{equation} 
\item Let $(1,x,0,z)$ $(x,z \in \bR \setminus \{0\})$ be the homogeneous coordinates of the point $P \in \SOL$. The parameters of the 
corresponding translation curve $t_{E_0P}$ are the following
\begin{equation}
\begin{gathered}
\phi=0~\text{or}~  \pi, ~\theta=\mathrm{arccot}\Big( \mp \frac{x}{(\mathrm{e}^{-z}-1)}\Big),\\
t=\frac{z}{\sin\theta}, ~ \text{where}  ~ -\pi/2\le \theta \le \pi/2, ~ t\in \bR^+.
\end{gathered} \tag{2.19}
\end{equation}
\item Let $(1,x,y,0)$ $(x,y \in \bR)$ be the homogeneous coordinates of the point $P \in \SOL$. The parameters of the 
corresponding translation curve $t_{E_0P}$ are the following
\begin{equation}
\begin{gathered}
\phi=\arccos\Big(\frac{x}{\sqrt{x^2+y^2}}\Big),~  \theta=0,\\
t=\sqrt{x^2+y^2}, ~ \text{where}  ~ -\pi < \phi \le \pi, ~ t\in \bR^+.
\end{gathered} \tag{2.20}
\end{equation}
\end{enumerate}
\end{lem}

The sphere of radius $R >0$ with centre at the origin, (denoted by $S^t_O(R)$), with the usual longitude and altitude parameters 
$\phi$ and $\theta$, can be obtained from (2.17) by replacing $t$ with $R$. It is easy to create the implicit equation of $S^t_{O}(R)$:
 \begin{equation}
\begin{gathered}
        \left(\frac{x z}{1-e^{-z}}\right)^2+\left(\frac{y z}{e^z-1}\right)^2+z^2=R^2
\end{gathered} \tag{2.21}
\end{equation}

\section{Geodesic triangles} \label{section3}
We consider $3$ points $A_1$, $A_2$, $A_3$ in the projective model of $\SOL$ space 
(see Section 2 and Fig.~2).
The {\it geodesic segments} $a_k$ between the points $A_i$ and $A_j$
$(i<j,~i,j,k \in \{1,2,3\}, k \ne i,j$) are called sides of the 
{\it geodesic triangle} with vertices $A_1$, $A_2$, $A_3$.
\begin{figure}[ht]
\centering
\includegraphics[width=12cm]{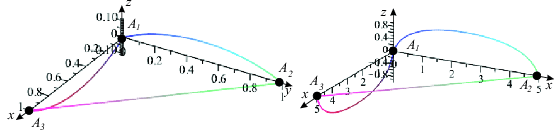}
\caption{Two different geodesic triangle: left with vertices $A_1=(1,0,0,0)$, $A_2=(1,1,0,0)$, $A_3=(1,0,1,0)$,
right with vertices $A_1=(1,0,0,0)$, $A_2=(1,5,0,0)$, $A_3=(1,0,5,0)$.}
\label{fig:Fig2}
\end{figure}
In Riemannian geometries the infinitesimal arc-length square (see (2.5)) and the corresponding metric tensor is used to define the angle $\theta$ between two geodesic curves.
If their tangent vectors in their common point are $\bu$ and $\bv$ and $g_{ij}$ are the components of the metric tensor then
\begin{equation}
\cos(\theta)=\frac{u_i g_{ij} v_j}{\sqrt{u_i g_{ij} u_j~ v_i g_{ij} v_j}} \tag{3.1}
\end{equation}
It is clear using (2.4) and the corresponding metric tensor, that
the angles are the same as the Euclidean ones at the origin.

We note here that the angle of two intersecting geodesic curves depend on the orientation of the tangent vectors. We will consider
the {\it interior angles} of the triangles that are denoted at the vertex $A_i$ by $\omega_i$ $(i\in\{1,2,3\})$.

\subsection{Horizontal-like isosceles triangles}
A geodesic triangle is called horizontal-like isosceles if its vertices lie on the 
$[x,y]$ plane with coordinates $A_1=(1,0,0,0)$, $A_2=(1,a,0,0)$, $A_3=(1,0,a,0)$ 
($a\in \mathbf{R}\setminus 0)$ (see Fig.~2). 
\begin{figure}[ht]
\centering
\includegraphics[width=6cm]{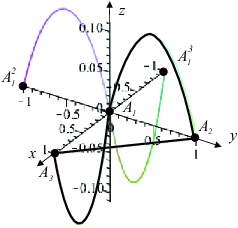}
\caption{The points $A_3$ and $A_1^3$ are antipodal related to the origin $E_0$ 
($A_1=(1,0,0,0)$,~$A_2=(1,0,0,\frac{1}{2})$,~$A_3=(1,4,0,\frac{1}{2})$.)}
\label{fig:Fig3}
\end{figure}
The geodesic segment $A_2A_3$ lies on a straight line, parallel to straight line $x+y=0$ (see (2.9)).  The geodesic segment $A_1A_2$ lies on the $[y,z]$ coordinate plane (see (2.8)
and the geodesic segment $A_2A_3$ lies on the coordinate plane $[x,z]$ as the image of geodesic segment $A_1A_2$ by the isometry given in (2.5) (2).

{\it In order to determine the further interior angles,
we define \emph{translations} $\bT_{A_i}$, $(i\in \{2,3\})$ as elements of the isometry group of $\SOL$, that
maps the origin $E_0=A_1$ onto $A_i$ (see Fig.~2).}
E.g. this $\SOL$ isometry (an Euclidean translation) $\bT_{A_3}$ and its inverse (up to a positive determinant factor) can be given by:
\begin{equation}
\bT_{A_3}=
\begin{pmatrix}
1 & a & 0 & 0 \\
0 & 1 & 0 & 0 \\
0 & 0 & 1 & 0 \\
0 & 0 & 0 & 1
\end{pmatrix} , ~ ~ ~
\bT_{A_3}^{-1}=
\begin{pmatrix}
1 & -a & 0 & 0 \\
0 & 1 & 0 & 0 \\
0 & 0 & 1 & 0 \\
0 & 0 & 0 & 1
\end{pmatrix} , \tag{3.2}
\end{equation}
and the images $\bT_{A_3}(A_i)$ of the vertices $A_i$ $(i \in \{1,2,3\})$ are the following:
\begin{equation}
\begin{gathered}
\bT^{-1}_{A_3}(A_1)=A_1^3=(1,-a,0,0),~\bT^{-1}_{A_3}(A_2)=A_2^3=(1,-a,a,0), \\ \bT^{-1}_{A_3}(A_3)=A_3^3=E_0=(1,0,0,0). \tag{3.3}
\end{gathered}
\end{equation}
Our aim is to determine angle sum $\sum_{i=1}^3(\omega_i)$ of the interior angles of the above right-angled
horizontal-like geodesic triangle $A_1A_2A_3$.
We mentioned that the angle of geodesic curves with common point at the origin $E_0$ is the same as the
Euclidean one. Therefore, it can be determined by usual Euclidean sense. 
Hence, $\omega_1$ is equal to the angle $\angle (g(E_0=A_1, A_3),g(E_0, A_2))$ (see Fig.~2)
where $g(A_1, A_3)$ and $g(A_1, A_2)$ are oriented geodesic curves.
Moreover, the translation $\bT_{A_3}$ is isometry
in $\SOL$ geometry thus
$\omega_3$ is equal to the angle $\angle (g(A_3^3, A_1^3),g(A_3^3, A_2^3)) $ (see Fig.~2.b)
where $g(A_3^3, A_1^3)$ and $g(A_3^3, A_2^3)$ are also oriented geodesic curves $(E_0=A_1=A_3^3)$.

We denote the oriented unit tangent vectors of the geodesic curves $g(E_0=A_1, A_i^j)$ with $\mathbf{t}_i^j$ where
$(i,j)\in \{(1,3),(2,3),(3,0),(2,0)\}$ and $A_3^0=A_3$, $A_2^0=A_2$.
The Euclidean coordinates of $\mathbf{t}_i^j$ (see Section 2) are :
\begin{equation}
\mathbf{t}_i^j=(\cos(\theta_i^j) \cos(\alpha_i^j), \cos(\theta_i^j) \sin(\alpha_i^j), \sin(\theta_i^j)). \tag{3.4}
\end{equation}
\begin{thm}
The sum $\Omega(a)=\sum_{i=1}^3(\omega_i(a))$ of the interior angles of a horizontal-like isosceles geodesic triangle $\pi < \Omega(a))< \frac{3\pi}{2}$ .
\end{thm}
\textbf{Proof:} If we consider the geodesic curve segment $g(A_1, A_2)$ then the corresponding parameters 
($\theta_2^0, \alpha_2^0, t_2^0$) follows from the formulas (2.8):
\begin{equation}
\theta_2^0=\arctan{\frac{a}{2}};~ \alpha_2^0=\frac{\pi}{2}; ~ t_2^0=\log \Big(\frac{a+\sqrt{a^2+4}}{\sqrt{a^2+4}-a} \Big). \tag{3.5}
\end{equation}
Therefore we obtain the corresponding tangent vector $\mathbf{t}_2^0$ and using the $\SOL$ isometry (2.5) (2) the tangent vector $\mathbf{t}_3^0$ 
and so their angle $\omega_1$.
\begin{equation}
\begin{gathered}
\mathbf{t}_2^0=(0, \cos(\arctan{\frac{a}{2}}), \sin(\arctan{\frac{a}{2}})),  \\
\mathbf{t}_3^0=(\cos(\arctan{\frac{a}{2}}), 0, -\sin(\arctan{\frac{a}{2}})).
\end{gathered} \notag
\end{equation}

It is clear, that the points $A_3^2$ and $A_2^3$ are antipodal to the origin $E_0=A_1$ (see Fig.~3) therefore, 
$\mathbf{t}_2^3=(1/\sqrt{2},-1/\sqrt{2},0)$ and $\mathbf{t}_3^2=(-1/\sqrt{2},1/\sqrt{2},0)$. 
Moreover, using the translations $\bT_{A_2}^{-1}$ and $\bT_{A_3}^{-1}$ (see (3.2)) we obtain that $\theta_1^2= \theta_2^0$ and $\theta_3^0=\theta_1^3$.  
From all this it follows that $\omega_2=\omega_3$ holds and
\begin{equation}
\begin{gathered}
\omega_1=\pi-\arccos \Big(\frac{4a}{(a^2+4)^{3/2}}\Big), ~ ~ \omega_2=\omega_3=\arccos \Big(\frac{2\sqrt{2}}{a^2+4} \Big).
\end{gathered} \tag{3.6}
\end{equation}

The interior angle sums $\Omega(a)=\sum_{i=1}^3(\omega_i(a))$ of horizontal-like isosceles geodesic triangles depend on 
parameter $a \in \mathbf{R}^+$. We investigated
this function and obtained the following results
$\Omega(a)$ ($a \in \mathbf{R}^+$) is a strictly increasing function and (see Fig.~4)
\begin{equation}
\begin{gathered}
\lim_{a \to 0}(\omega_1(a))=\frac{\pi}{2}, ~ \lim_{a \to 0}(\omega_2(a))=\lim_{a \to 0}(\omega_3(a))=\frac{\pi}{4},\\
\lim_{a \to \infty}(\omega_1(a))= \lim_{a \to \infty}(\omega_2(a))=\lim_{a \to \infty}(\omega_3(a))=\frac{\pi}{2}, \\
\lim_{a \to 0}(\Omega(a))=\pi, ~ \lim_{a \to \infty}(\Omega(a))=\frac{3\pi}{2}.
\end{gathered} \notag 
\end{equation} \ \  \ \ \ $\square$
\begin{conj}
The sum of the interior angles of any horizontal-like geodesic triangle is greater than $\pi$.
\end{conj}
\begin{figure}[ht]
\centering
\includegraphics[width=6cm]{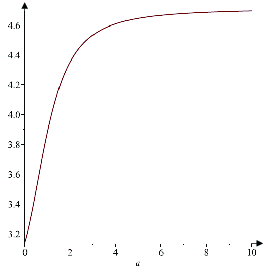}
\caption{}
\label{fig:Fig4}
\end{figure}
In the following table we summarize some numerical data of geodesic triangles for given parameters:
\medbreak
\centerline{\vbox{
\halign{\strut\vrule\quad \hfil $#$ \hfil\quad\vrule
& \quad \hfil $#$ \hfil\quad\vrule &\quad \hfil $#$ \hfil\quad\vrule &\quad \hfil $#$ \hfil\quad\vrule \cr
\noalign{\hrule}
\multispan4{\strut\vrule\hfill\bf Table 1, ~ Horizontal-like isosceles geodesic triangles \hfill\vrule} \cr
\noalign{\hrule}
\noalign{\vskip2pt}
\noalign{\hrule}
a & \omega_1 & \omega_2(a)=\omega_3(a)  & \Omega(a)=\sum_{i=1}^3(\omega_i(a))  \cr
\noalign{\hrule}
\rightarrow 0 & \rightarrow \pi/2 & \rightarrow \pi/4 & \rightarrow \pi \cr
\noalign{\hrule}
1 & 1.93668 &  0.96953 & 3.87574 \cr
\noalign{\hrule}
5 & 1.69922 &  1.47311 & 4.64543 \cr
\noalign{\hrule}
50 & 1.57239 & 1.56967 & 4.71173 \cr
\noalign{\hrule}
1000 & 1.57080 & 1.57079 & 4.71239 \cr
\noalign{\hrule}
\rightarrow \infty & \rightarrow \pi/2 & \rightarrow \pi/2 & \rightarrow 3\pi/2 \cr
\noalign{\hrule}
}}}
\medbreak
\subsection{Hyperbolic-like geodesic triangles}
A geodesic triangle is hyperbolic-like if its vertices lie in the $[y,z]$ or $[x,z]$ coordinate plane of the model or a plane parallel with it. 
In this section we analyse the interior angle sums of the hyperbolic-like geodesic triangles. 

It is clear by the Remark 2.2 that the $[x,z]$ plane (and so the plane $[y,z]$) is a convexly embedded copy 
of $\mathbf{H}^2$ in $\SOL$ and similar statement is true for the planes $x=c$, $y=c$ $(c\in \mathbf{R})$. Therefore we can formulate the following
\begin{thm}
The interior angle sum $\Omega=\sum_{i=1}^3(\omega_i)$ of a hyperbolic-like geodesic triangle is less than $\pi$.  \ \ \ $\square$
\end{thm}
For example, we study the hyperbolic-like geodesic triangle with wertices $A_1=(1,0,0,0),~A_2=(1,0,a,0),~A_3=(1,0,0,a)$ ($a \in \mathbf{R}^+$) (see Fig.~5).
\begin{figure}[ht]
\centering
\includegraphics[width=7cm]{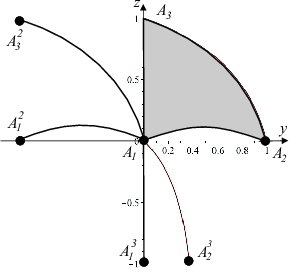}
\caption{The hyperbolic-like geodesic triangle with vertices $A_1=(1,0,0,0),~A_2=(1,0,1,0),~A_3=(1,0,0,1)$ and
the its translated geodesic edges $A_1A_1^2$, $A_1A_1^3$, $A_1A_2^3$ and $A_1A_3^2$.}
\label{fig:Fig5}
\end{figure}
The geodesic segment $A_1A_3$ lies on the $z$ axis (see (2.10)) and the geodesic segments $A_1A_2$, $A_2A_3$ lie in the $[y,z]$ coordinate plane (see (2.8)). 

In order to determine the interior angles of hyperbolic-like geodesic triangle $A_1A_2A_3$ 
similarly to the horizontal-like case we use the Euclidean-like translation $\bT_{A_2}$ (see (3.2)) and define the translation $\bT_{A_3}$ that 
maps the origin $E_0=A_1$ onto $A_3$ that can be given by:
\begin{equation}
\bT_{A_3}=
\begin{pmatrix}
1 & 0 & 0 & a \\
0 & e^{-a} & 0 & 0 \\
0 & 0 & e^{a} & 0 \\
0 & 0 & 0 & 1
\end{pmatrix} , ~ ~ ~
\bT_{A_3}^{-1}=
\begin{pmatrix}
1 & 0 & 0 & -a \\
0 & e^{a} & 0 & 0 \\
0 & 0 & e^{-a} & 0 \\
0 & 0 & 0 & 1
\end{pmatrix}. \tag{3.7}
\end{equation}
We obtain that the images $\bT^{-1}_{A_3}(A_i)$ of the vertex $A_i$ $(i \in \{1,2,3\}),$
are the following (see also Fig.~5):
\begin{equation}
\begin{gathered}
\bT^{-1}_{A_3}(A_1)=A_1^3=(1,0,0,-a),~\bT^{-1}_{A_3}(A_2)=A_2^3=(1,0,a \cdot e^{-a},-a), \\ \bT^{-1}_{A_3}(A_3)=A_3^3=A_1=E_0=(1,0,0,0). \tag{3.8}
\end{gathered}
\end{equation}
We study similarly to the above horizontal-like case the sum $\Omega(a)=\sum_{i=1}^3(\omega_i(a))$ of the interior angles of the above hyperbolic-like geodesic triangle $A_1A_2A_3$
The computations are similar to the former case, I will not detail here.
\begin{equation}
\begin{gathered}
\theta_2^0=\arctan{\frac{a}{2}};~ \alpha_2^0=\frac{\pi}{2}; ~ t_2^0=\log \Big(\frac{a+\sqrt{a^2+4}}{\sqrt{a^2+4}-a} \Big),\\ \tag{3.9}
\theta_3^2=\arctan{ \frac{e^{2a}+a^2-1}{2a}}; ~ \alpha_3^2=\frac{\pi}{2}; ~ 
t_2^3=t_2^3=-\log{2}+a+\\ +\log \Big( \frac{(a^2+1)^2 e^{-2a}+2a^2+e^{2a}-2}{\sqrt{a^4+2a^2(e^{2a}+1)+e^{4a}-2e^{2a}+1}}+e^{-2a}(a^2+1)+1 \Big);\\
\theta_2^3=\arctan{ \frac{e^{2a}-a^2-1}{-2ae^a}}; ~ \alpha_2^3=\frac{\pi}{2}; ~ 
\end{gathered} 
\end{equation}
The interior angle sums $\Omega(a)=\sum_{i=1}^3(\omega_i(a))$ of the above hyperbolic-like geodesic triangles depend on 
parameter $a \in \mathbf{R}^+$ and investigating this function (see (3.10)) we obtain the following results
\begin{lem}
The sum $\Omega(a)=\sum_{i=1}^3(\omega_i(a))$ of the interior angles of hyperbolic-like geodesic triangles with vertices $A_1=(1,0,0,0),~A_2=(1,0,a,0),~A_3=(1,0,0,a)$
$0 < \Omega(a))< \pi$ .
\begin{equation}
\begin{gathered}
\Omega(a)=\arctan{ \frac{e^{2a}+a^2-1}{2a}}-2\arctan{\frac{a}{2}}-\arctan{ \frac{e^{2a}-a^2-1}{-2ae^a}}+\pi
\end{gathered} \notag
\end{equation}
$\Omega(a)$ ($a \in \mathbf{R}^+$) is a strictly decreasing function and (see Fig.~6)
\begin{equation}
\begin{gathered}
\lim_{a \to 0}(\omega_1(a))=\frac{\pi}{2}, ~ \lim_{a \to 0}(\omega_2(a))=\frac{\pi}{4}; \lim_{a \to 0}(\omega_3(a))=\frac{\pi}{4},\\
\lim_{a \to \infty}(\omega_1(a))= 0; ~ \lim_{a \to \infty}(\omega_2(a))=0; ~ \lim_{a \to \infty}(\omega_3(a))=0, \\
\lim_{a \to 0}(\Omega(a))=\pi, ~ \lim_{a to \infty}(\Omega(a))=0.
\end{gathered} \notag 
\end{equation} \ \  \ \ \ $\square$
\end{lem}

\begin{figure}[ht]
\centering
\includegraphics[width=6cm]{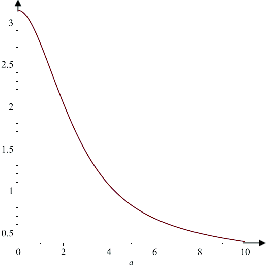}
\caption{}
\label{}
\end{figure}

We can determine the interior angle sum of arbitrary hyperbolic-like geodesic triangle similarly to the fibre-like case.
In the following table we summarize some numerical data of geodesic triangles for given parameters:
\medbreak
\centerline{\vbox{
\halign{\strut\vrule\quad \hfil $#$ \hfil\quad\vrule
&\quad \hfil $#$ \hfil\quad\vrule &\quad \hfil $#$ \hfil\quad\vrule &\quad \hfil $#$ \hfil\quad\vrule &\quad \hfil $#$ \hfil\quad\vrule
\cr
\noalign{\hrule}
\multispan5{\strut\vrule\hfill\bf Table 2, ~ ~ Hyperbolic-like geodesic triangles \hfill\vrule}%
\cr
\noalign{\hrule}
\noalign{\vskip2pt}
\noalign{\hrule}
a & \omega_1 & \omega_2(a) & \omega_3(a) & \Omega(a)=\sum_{i=1}^3(\omega_i(a))  \cr
\noalign{\hrule}
\rightarrow 0 & \rightarrow \pi/2 & \rightarrow \pi/4 & \rightarrow \pi/4 & \rightarrow \pi \cr
\noalign{\hrule}
1 & 1.10715 & 0.84281 & 0.78979 & 2.73975 \cr
\noalign{\hrule}
5 & 0.38051 &  0.38005 &  0.06736 & 0.82792 \cr
\noalign{\hrule}
50 & 0.03998 &  0.03998 & 0.00000 & 0.07996 \cr
\noalign{\hrule}
1000 & 0.00200 &  0.00200 & 0.00000 & 0.00400 \cr
\noalign{\hrule}
\rightarrow \infty &  \rightarrow 0 & \rightarrow 0 & \rightarrow 0 & \rightarrow 0 \cr
\noalign{\hrule}
}}}
\medbreak
\subsection{Geodesic triangles with interior angle sum $\pi$}
In the above sections we discuss the horizontal- and hyperbolic-like geodesic triangles and proved that there are 
geodesic triangles whose angle sum $\Omega(a)=\sum_{i=1}^{3}(\omega_i(a))$ is greater or less than $\pi$. $\Omega(a)=\sum_{i=1}^{3}(\omega_i(a))=\pi$ is
realized if the parameter $a$ tends to $0$. We prove the following
\begin{lem}
There is geodesic triangle $A_1A_2A_3$ with interior angle sum $\pi$ where its vertices are {\it proper}
(i.e. $A_i \in \SOL$, vertices are not in infinity ($i \in \{1,2,3\}$) or not tend to the origin).
\end{lem}
{\bf{Proof:}} We consider a horizontal-like isosceles geodesic triangle with vertices $A_1=E_0=(1,0,0,0)$, $A_2^h=(1,0,a,0)$
$A_3^h=(1,a,0,0)$ and a hyperbolic-like geodesic triangle with vertices
$A_1=E_0=(1,0,0,0)$, $A_2^v=(1,0,0,a)$ $A_3^v=(1,0,0,a)$
($a \in \mathbf{R}^+$) . We cosider the straight line segments (in Euclidean sense) $A_3^hA_3^v$, ($A_3^h A_3^f \subset \SOL$).

We consider a geodesic triangle $A_1A_2A_3(t)$ where $A_3(t) \in A_3^h A_3^v$ $(t\in [0,1])$. $A_3(t)$ are moving on the
segments $A_3^hA_3^v$. If $t=0$ then $A_3(0)=A_3^h$ and if $t=1$ then $A_3(1)=A_3^v$.

Similarly to the above cases the internal angles of the geodesic triangle $A_1A_2(t)$ $A_3(t)$ are denoted by $\omega_i(t)$ $(i \in \{1,2,3\})$.
The angle sum $\sum_{i=1}^{3}(\omega_i(0)) > \pi$ and $\sum_{i=1}^{3}(\omega_i(1)) < \pi$ (see Lemma 3.1 and Lemma 3.3-4) . Moreover the angles $\omega_i(t)$ change
continuously if the parameter $t$ run in the interval $[0,1]$. Therefore there is a $t_E\in (0,1)$ where $\sum_{i=1}^{3}(\omega_i(t_E)) = \pi$.
~ ~ $\square$

We obtain by the Lemmas of this Section the following
\begin{thm}
The sum of the interior angles of a geodesic triangle of $\SOL$ space can be larger than, less than or equal to $\pi$.
\end{thm}
\section{Translation-like isoptic surfaces in $\SOL$}

\begin{defn} The $\SOL$ translation-like $\alpha$ -- isoptic surface of a translation-like segment $\overline{A_1A_2}$ is the locus of points $P$ 
for which the internal angle at $P$ in the translation-like triangle, formed by $A_1,$ $A_2$ and $P$ is $\alpha.$ If $\alpha$ is the right angle, 
then it is called the translation-like Thaloid of $\overline{A_1A_2}.$
\label{defisoszakaszN}
\end{defn}

In the rest of this study, we will focus on the isoptic surface of the 
translation-like segment in $\SOL$ geometry, which in the projective model is far from straight, but a parametric curve described in (2.17). 

We emphasize here that the section itself does not appear in our calculations, 
we only deal with the endpoints. 
An interesting question beyond this study is how the ruled surface, {or more precisely in this case, how the triangular surface}  
looks like generated by the curves drawn from 
the outer point to all points 
of the section. Thus the angle can really be considered planar {in $\SOL$ sense} or any 
non desirable intersection occurs between the segment and the rays. 
The section itself and the rays can be translation-like or geodesic-like as well.
Some of these questions arise in \cite{Sz22}.

Let us be given a point in $\SOL$ with its homogeneous coordinates $A_1(1,a,b,c)$. We would like to imagine or segment symmetric respected to the origin, therefore we consider $A_2$ as a point on the translation curve from the origin to $A_1$ antipodal. So that we replace $t$ with $-t$ in (2.17) to obtain the homogeneous coordinates of $A_2$. 

\begin{equation}
\begin{gathered}
        \left\{ \begin{array}{ll} 
        x(\phi,\theta,-t)=-\cot{\theta} \cos{\phi} (e^{t \sin{\theta}}-1)= - x(\phi,\theta,t) e^{z(\phi,\theta,t)}, \\
        y(\phi,\theta,-t)=\cot{\theta} \sin{\phi} (e^{-t \sin{\theta}}-1)= 
        - y(\phi,\theta,t) e^{-z(\phi,\theta,t)}, \\
        z(\phi,\theta,-t)=-t \sin{\theta}= - z(\phi,\theta,t).
        \end{array} \right. \\
        \text{If $\theta=0$ then}: ~  x(\phi,\theta,-t)=-t\cos{\phi}=-x(\phi,\theta,t) , \\ y(\phi,\theta,-t)=-t \sin{\phi}=-y(\phi,\theta,t),  \ z(\phi,\theta,-t)=0.
\end{gathered}
\end{equation}

Now that $A_2=(1,-a e^c,-be^{-c},-c)$ is known, considering a point $P(1,x,y,z)$, we can determine the angle $A_1PA_2\angle$ along the 
procedure described in the previous section. We apply $\bT_{P}^{-1}$ to all three points. 
This transformation preserves the angle $A_1PA_2\angle$ and pulls
back $P$ to the origin, 
hence the angle in question seems in real size. 
We get $\bT_{P}^{-1}$ by inverting the matrix in (2.3):

\begin{equation}
\bT_{P}^{-1}=
\begin{pmatrix}
1 & -x e^z & -y e^{-z} & -z \\
0 & e^z & 0 & 0 \\
0 & 0 & e^{-z} & 0 \\
0 & 0 & 0 & 1 
\end{pmatrix},
\end{equation}

\begin{equation}
\begin{gathered}
\bT^{-1}_{P}(P)=(1,0,0,0), \\
\bT^{-1}_{P}(A_1)=(1,e^z(a-x),e^{-z}(b-y),c-z),\\
\bT^{-1}_{P}(A_2)=(1,-e^z(ae^c+x),-e^{-z}(be^{-c}+y),-(z+c)).
\label{eltoltak}
\end{gathered}
\end{equation}

According to (2.16) and (2.17), the tangent of the translation curve between the origin and a point $P(1,x,y,z)$ at the origin can be obtained by the following formulas:

\begin{equation}
\bt=(u,v,w)=\left(\dfrac{x z}{1-e^{-z}},\dfrac{y z}{e^{z}-1},z\right)
\label{erintfor}
\end{equation}

Let us denote with $\bt_1$ and $\bt_2$ the tangents of the translation curves to $\bT_{P}^{-1}(A_1)$ and $\bT_{P}^{-1}(A_2)$ from the origin $E_0=\bT_{P}^{-1}(P)$ at the origin. We can calculate these tangents by applying \ref{erintfor} to \ref{eltoltak}.

\begin{equation}
\begin{gathered}
\bt_{1}=-\dfrac{z-c}{e^{\frac{z-c}{2}}-e^{-\frac{z-c}{2}}}
\left((x-a)e^{\frac{z+c}{2}},(y-b)e^{-\frac{z-c}{2}},
e^{\frac{z-c}{2}}-e^{-\frac{z-c}{2}}\right)\\
\bt_{2}=-\dfrac{z+c}{e^{\frac{z+c}{2}}-e^{-\frac{z+c}{2}}}
\left((x+ae^c)e^{\frac{z-c}{2}},(y+b e^c)e^{-\frac{z-c}{2}},
e^{\frac{z+c}{2}}-e^{-\frac{z+c}{2}}\right)
\end{gathered}
\end{equation}

Finally, fixing the angle of the $\bt_1$ and $\bt_2$ to $\alpha,$ we get the {\it translation-like} $\alpha$ -- isoptic surface of $\overline{A_1A_2}.$ 

\begin{thm}
Given a {\it translation-like} segment in the $\SOL$ geometry by its endpoints $A_1=(1,a,b,c)$ and $A_2=(1,-ae^c,-be^{-c},-c).$ Then the translation-like 
$\alpha$ -- isoptic surface of the translation-like segment $\overline{A_1A_2}$ have the implicit equation:
{\footnotesize
\begin{equation}
\begin{gathered}
\cos(\alpha)=\displaystyle {\frac{\bt_1\cdot\bt_2}{|\bt_1|\cdot|\bt_2|}}=\\
\displaystyle{\frac {e^z\left(x+a e^c\right)\left(x-a \right) + e^{-z}\left(y+b e^{-c}\right)\left( y-b \right)+ 4\sinh{\frac{z+c}{2}}\sinh{\frac{z-c}{2}} }
{\splitfrac{ \left( e^{z+c}(x-a)^2+e^{-(z+c)}(y-b)^2+ 4\sinh^2{\frac{z-c}{2}}\right)^\frac{1}{2}\cdot}{\cdot\left( e^{z-c}\left(x+a e^c\right)^2+e^{-(z-c)}\left(y+be^{-c}\right)^2+ 4\sinh^2{\frac{z+c}{2}}\right)^{\frac{1}{2}}\ \ \square}
}}
\label{eq:soliso}
\end{gathered}
\end{equation}}
\end{thm}

\begin{figure}
\centering
\includegraphics[width=6cm]{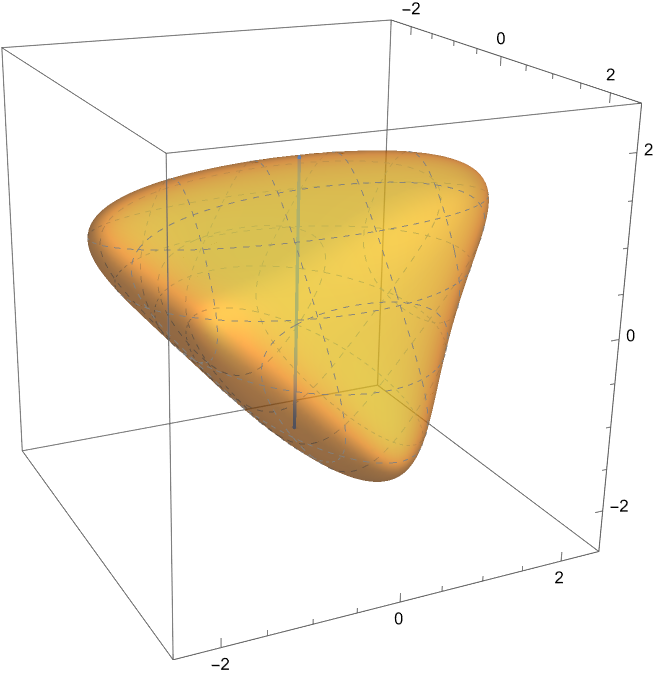}
\includegraphics[width=6cm]{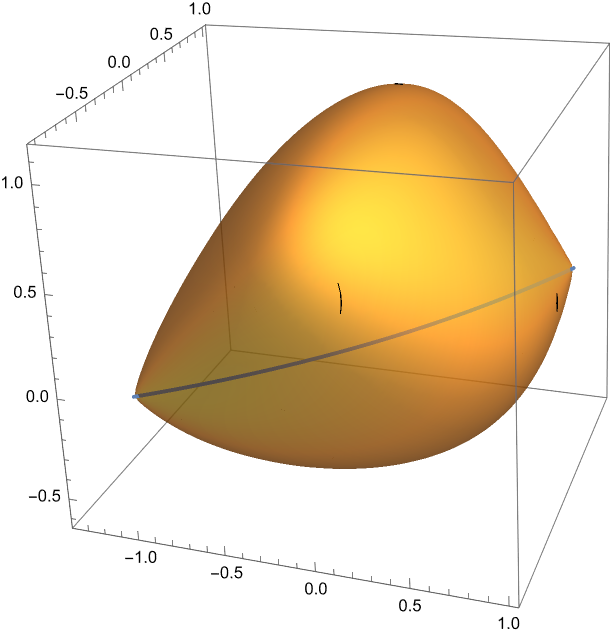}
\caption{Isoptic surface with $A_1=(0,0,1.6)$ and $\alpha=\frac{\pi}{2}$ (left) and with $A_1=(1,0.3,0.4)$ and $\frac{2\pi}{3}$ (right).}
\label{fig:Fig7}
\end{figure}

\begin{figure}
\centering
\includegraphics[width=6cm]{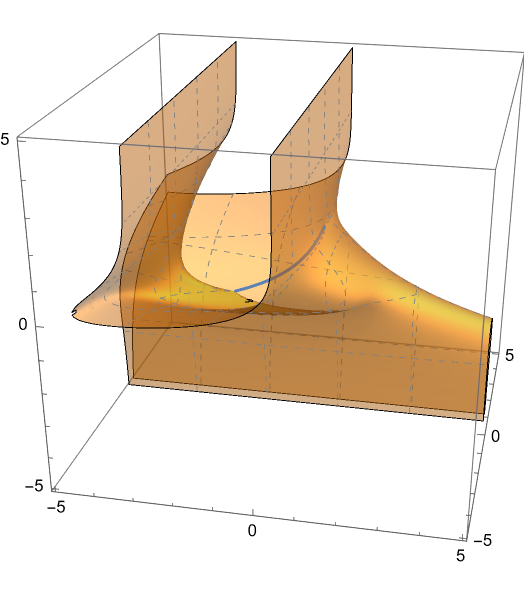}
\includegraphics[width=6cm]{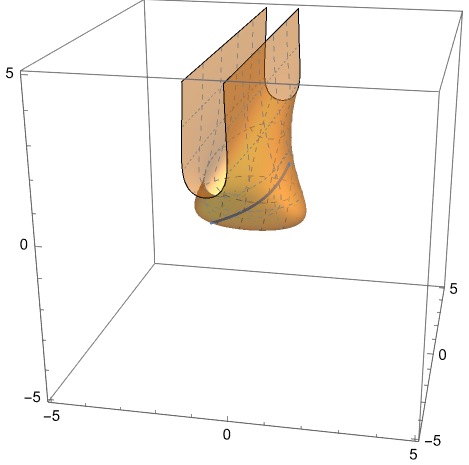}
\caption{Isoptic surface with $A_1=(0.6,0.4,1)$ and $\alpha=\frac{\pi}{3}$ (left) and $\frac{\pi}{2}$ (right).}
\label{fig:Fig8}
\end{figure}
On Fig.~\ref{fig:Fig7}, one can see some isoptic surfaces to a general translation-like segment in $\SOL$ geometry. The right side shows the isoptic surface to an obtuse angle, 
the left side shows the translation-like Thaloid of another segment.

In the following, we will examine some special properties of the above equation. As we have seen e.g. in \cite{Cs23,Cs-Sz23}, it is expected that the Thales surface $(\alpha=\frac{\pi}{2})$ of a specially positioned section  coincides with the origin-centered sphere of the given geometry. This property was fulfilled in all other Thurston geometries in the translational case. In light of the above, the following lemma is quite surprising

\begin{lem}
Given a {\it translation-like} segment in the $\SOL$ geometry by its endpoints  $A_1=(1,a,b,c)$ and $A_2=(1,-ae^c,-be^{-c},-c).$  Then the \it{translation-like Thaloid} 
of this line segment will be no translation sphere for any values of $a,$ $b$ and $c.$
\label{lemiso}
\end{lem}
\textbf{Proof:} At first, it can be easy seen that (2.21) is equivalent with the following implicit equation:
\begin{equation}
\left(\frac{\frac{z}{2}}{\sinh{\frac{z}{2}}}\right)^2\left(e^z x^2+e^{-z} y^2+4\sinh^2{\frac{z}{2}}\right)-R^2=0,
\label{eq:sols}
\end{equation}
where $R$ is the translational distance of $A_1$ and the origin can be computed from $a,$ $b$ and $c.$

For the equation of the Thaloid, we can make the numerator of the left side in (\ref{eq:soliso}) equal to 0:
\begin{equation}
 e^z\left(x+a e^c\right)\left(x-a \right) + e^{-z}\left(y+b e^{-c}\right)\left( y-b \right)+ 4\sinh^2{\frac{z}{2}}-\sinh^2{\frac{c}{2}}=0   
 \label{eq:thal}
\end{equation}
Now, we can compare (\ref{eq:sols}) and (\ref{eq:thal}). It is easy to see that in (\ref{eq:thal}), the coefficient of $x e^z$ is 0 like in (\ref{eq:sols}) if $a(e^c-1)=0.$ Similarly, the coefficient of $y e^{-z}$ is 0, if $b(e^{-c}-1)=0.$ We can distinguish 2 cases than:
\begin{enumerate}
    \item[1. case] $c=0$ \\
    According to (2.20), the radius of the sphere is $R=\sqrt{a^2+b^2}.$ Replacing $c$ with 0 in (\ref{eq:sols}) and in (\ref{eq:thal}), we get the following equations:
    \begin{equation}
 \mathrm{Thaloid:}\ \ e^z x^2+ e^{-z}y^2+ 4\sinh^2{\frac{z}{2}}=e^z a^2+e^{-z}b^2 
\end{equation}
\begin{equation}
 \mathrm{Sphere:}\ \ e^z x^2+ e^{-z}y^2+ 4\sinh^2{\frac{z}{2}}=(a^2+b^2 )\frac{4\sinh^2{\frac{z}{2}}}{z^2}
\end{equation}
The above equations were rearranged so that their left sides are equal. It is easy to see that the right sides cannot be equal for all $z.$
 \item[2. case] $a=0\wedge b=0$ \\
  According to (2.19), the radius of the sphere is $R=c.$ Replacing $a$ and $b$ with 0 in (\ref{eq:sols}) and in (\ref{eq:thal}), we get the following equations:
    \begin{equation}
 \mathrm{Thaloid:}\ \ e^z x^2+ e^{-z}y^2+ 4\sinh^2{\frac{z}{2}}=4\sinh^2{\frac{c}{2}}
\end{equation}
\begin{equation}
 \mathrm{Sphere:}\ \ e^z x^2+ e^{-z}y^2+ 4\sinh^2{\frac{z}{2}}=c^2\frac{4\sinh^2{\frac{z}{2}}}{z^2}
\end{equation}
The above equations were rearranged so that their left sides are equal. It is easy to see that the right sides cannot be equal for all $z.$ \ \ \ $\square$
\end{enumerate} 

\begin{rem}
We can also prove this lemma by using the invariancy of the sphere  under the isometries, described in (2.5). The Thaloid will be invaraint for $y\leftrightarrow -y$ transformation if and only if $b=0,$ or $c=0.$ We can evaluate the invariancy of other isometry $x\leftrightarrow y,$ $z\leftrightarrow-z$ regarding these cases. The Thaloid is invariant under this isometry for $b=0$, if and only if $a=0$ and for $c=0$, if and only if $a=\pm b.$
\end{rem}

It is also interesting based on Figure $\ref{fig:Fig8}$ that, in some cases, the isoptic surface is not a closed surface. That is why it may be interesting for us to examine the limit of (4.6), when $z$ tends to $\infty$ or $-\infty:$
\begin{equation}
    \lim_{z\rightarrow\infty}\cos(\alpha)=\frac{x^2+a \left(e^c-1\right) x+1-a^2 e^c}{\sqrt{e^{2 c} (a-x)^2+(a-x)^2 \left(a e^c+x\right)^2+e^{-2 c} \left(a
   e^c+x\right)^2+1}}
    \label{eq:liminfp}
\end{equation}
\begin{equation}
    \lim_{z\rightarrow-\infty}\cos(\alpha)=\frac{y^2+b \left(e^{-c}-1\right) y+1-b^2 e^{-c}}{\sqrt{e^{-2 c} (b-y)^2+(b-y)^2 \left(b e^{-c}+y\right)^2+e^{2 c} \left(b
   e^{-c}+y\right)^2+1}}
       \label{eq:liminfm}
\end{equation}

Since the two limits are very similar, we will only deal with the $z\rightarrow\infty$ case and leave the discussion of the other to the gentlest reader.

In such a positive fiber-like direction, the isoptic surface cannot have a point if the calculated limit value is less than the cosine of the given $\alpha$ angle. Therefore it is crucial, to find the minima of (\ref{eq:liminfp}). Simple calculus can be used to determine the $x_{1,2,3}$ location of the extreme values of (\ref{eq:liminfp}). If $c=0$ then $x_1=0,$ and the only minimum is $\dfrac{1-a^2}{1+a^2}$ otherwise

\begin{equation}
    x_{2}=\frac{a e^c}{e^c-1} \ \mathrm{and} \ x_{1,3}=\frac{a e^c\pm\sqrt{a^2 e^{2 c}+\left(e^c-1\right)^2}}{e^c-1}
    \label{eq:extrema}
\end{equation}

It is easy to check that the value of the right side in (\ref{eq:liminfp}) at $x_2$ is $1$, which is definitely a maximum. Furthermore, the limit of it taken at $\infty$ and $-\infty$ are both $1.$ Its continuity grants that $x_1$ and $x_3$ are minimum places. Exactly which minimum will be the absolute minimum of the function requires further discussion. We omit this calculation for reasons of scope and only support the following statement with a figure. 
\begin{lem}
Given a {\it translation-like} segment in the $\SOL$ geometry by its endpoints $A_1=(1,a,b,c)$ and $A_2=(1,-ae^c,-be^{-c},-c).$ Then the translation-like 
$\alpha$ -- isoptic surface of the translation-like segment $\overline{A_1A_2}$ is an infinite surface, if
{\footnotesize
\begin{equation}
\begin{gathered}
\cos(\alpha)>\displaystyle {\frac{e^c}{\sqrt{\left(a^2+1\right) e^{2 c}-2 e^c+1}}\cdot}\\
\displaystyle{\frac{ 2-4 e^c\pm a \sqrt{\left(a^2+1\right) e^{2 c}-2 e^c+1}+e^{2 c} \left(2
   a^2+2\pm a \sqrt{\left(a^2+1\right) e^{2 c}-2 e^c+1}\right)}{\splitfrac{\left(\left(a^2+1\right) e^{6 c}+e^{2 c} 
   \left(a^2+3\pm4 a \sqrt{\left(a^2+1\right) e^{2 c}-2 e^c+1}\right)\right.+}{\left.+e^{4 c} \left(6 a^2+3\pm4 a \sqrt{\left(a^2+1\right) e^{2 c}-2 e^c+1}\right)-2 e^c-4
   e^{3 c}-2 e^{5 c}+1\right)^\frac{1}{2},}}}
\label{eq:soliso2}
\end{gathered}
\end{equation}
where $\pm$ is $-$ if $a>0$ and $+$ if $a<0.$ $\square$}
\end{lem}

\begin{figure}
    \centering
    \includegraphics[width=8cm]{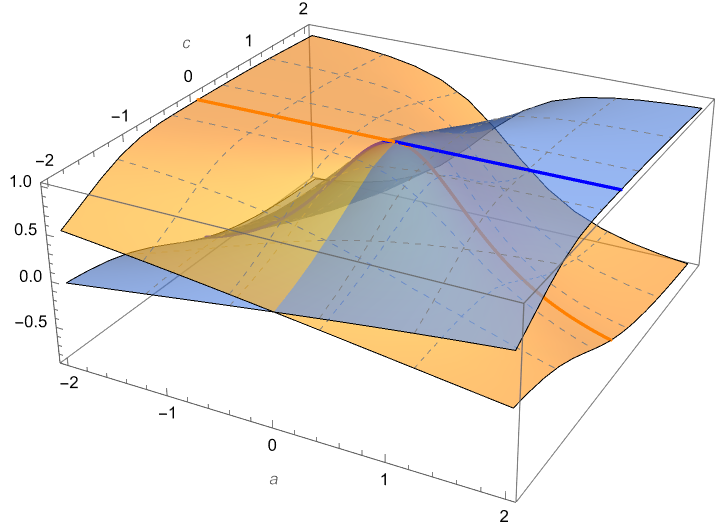}
    \caption{Right side of the inequality in (\ref{eq:soliso2}). Orange: $\pm$ is $-$ Blue: $\pm$ is $+$}
    \label{fig:limit}
\end{figure}

\begin{rem}
    Similar lemma can be given for $b,$ examining the $z\rightarrow-\infty$ case. Furthermore, it may be interesting to examine the limit values in other directions besides the $z$ axis, but our aim with the above examination was to prove that the isoptic surface is not always closed and not to perform a whole discussion.
\end{rem}

\section*{Data Availability Statement}
Data sharing not applicable to this article as no databases were generated or analyzed during the current study.

\section*{Conflict Of Interest Statement}
The authors have no conflict of interest to declare that are relevant to the content of this study.


\begin{thebibliography}{12}
%
\bibitem{BT}
{Brieussel,~J.~--~Tanaka,~R.,}
Discrete random walks on the group $\SOL$.
\textit{Isr. J. Math.,} {\bf 208/1}, 291-321 (2015).
%
\bibitem{B}
{Brodaczewska,~K.,}
Elementargeometrie in $\NIL$.
\textit{Dissertation (Dr. rer. nat.) Fakult\"at Mathematik und Naturwissenschaften der Technischen Universit\"at Dresden}
(2014).
%
\bibitem{CaMoSpSz}
{Cavichioli,~A.~--~Moln\'ar,~E.~--~Spaggiari,~F.~--~Szirmai,~J.,}
Some tetrahedron manifolds with $\SOL$ geometry.
\textit{J. Geom.,} {\bf 105/3}, 601-614 (2014).
%
\bibitem{harom} Cie\'slak, W.~--~Miernowski, A.~--~ Mozgawa, W.: 
Isoptics of a Closed Strictly Convex Curve, 
\textit{Lect. Notes in Math.}, {\bf 1481} (1991), pp. 28--35. 
%
\bibitem{tizenketto} Cie\'slak, W.~--~Miernowski, A.~--~Mozgawa, W.:
Isoptics of a Closed Strictly Convex Curve II, 
\textit{Rend. Semin. Mat. Univ. Padova} \textbf{96} (1996), 37--49.
%
\bibitem{Ch}
Chavel,~I.,
{\it Riemannian Geometry: A Modern Introduction}.  
Cambridge Studies in Advances Mathematics, (2006).
%
\bibitem{CsSz16}
Csima,~G.~--~Szirmai,~J.,
Interior angle sum of translation and geodesic triangles in $\SLR$ space. 
{\it Filomat}, {\bf 32/14}, 5023--5036, (2018). 
%
%
\bibitem{Cs23}
{{Csima,~G.:}
Isoptic surfaces of segments in $\SXR$ and
$\HXR$ geometries.
\emph{J. Geom.}
(2024), doi: 10.1007/s00022-023-00699-x.}
%
\bibitem{Cs-Sz13}
{{Csima,~G.~--~Szirmai,~J.:}
On the isoptic hypersurfaces in the n-dimensional Euclidean space.
\emph{KoG (Scientific and professional journal of Croatian Society for Geometry and Graphics)}
{\bf 17} (2013) 53--57.}
%
\bibitem{Cs-Sz14}
{{Csima,~G.~--~Szirmai,~J.:}
Isoptic curves of conic sections in constant curvature geometries.
\emph{Mathematical Communications}
{\bf 19} (2014) 277--290.}
%
\bibitem{Cs-Sz16}
{{Csima,~G.~--~Szirmai,~J.:}
Isoptic surfaces of convex polyhedra.
\emph{Computer Aided Geometric Design (CAGD)}
(2016) (to appear), DOI: 10.1016/j.cagd.2016.03.001.}
%
\bibitem{Cs-Sz16-1}
{Csima,~G.~--~Szirmai,~J.:}
Isoptic curves of generalized conic sections in the hyperbolic plane.
\emph{Ukrainian Mathematical Journal}, {\bf 71/12} (2020), 1929-1944, doi: 10.1007/s11253-020-01756-3.
%
\bibitem{Cs-Sz23}
{Csima,~G.~--~Szirmai,~J.:}
Translation-like isoptic surfaces and angle sums of translation triangles in $\NIL$ geometry.
\emph{Results Math.}, (2023) 78:194, DOI: 10.1007/s00025-023-01961-z,  arXiv: 2302.07653.
%
\bibitem{H} Holzm\"uller,~G.: Einf\"uhrung in die Theorie der isogonalen Verwandtschaft,
\textit{B.G. Teuber}, Leipzig-Berlin, (1882).
%
\bibitem{KN}
Kobayashi,~S.~--~Nomizu,~K.,
{\it Fundation of differential geometry, I.}.  Interscience, Wiley, New York (1963).
%
\bibitem{KV}
{Kotowski,~M.~--~Vir\'ag,~B.,}
Dyson's spike for random Schroedinger operators and Novikov-Shubin invariants of groups.
\textit{Manuscript (2016)} arXiv:1602.06626.
%
\bibitem{Kunkli} Kunkli, R.~--~Papp, I.~--~Hoffmann, M.: 
Isoptics of Bézier curves, 
\textit{Comput. Aided Geom. Design} \textbf{30} (2013), 78--84.
%
\bibitem{Kurusa} Kurusa, \'A.: 
Is a convex plane body determined by an isoptic?,  
\textit{Beitr. Algebra Geom.} {\bf 53} (2012), 281--294.
%
\bibitem{Loria} Loria, G.:Spezielle algebrische und transzendente ebene Kurven. 1 \& 2,
B.G. Teubner, Leipzig-Berlin, (1911).
%
\bibitem{nyolc} Michalska, M.: 
A sufficient condition for the convexity of the area of an isoptic curve of an oval, 
\textit{Rend. Semin. Mat. Univ. Padova} \textbf{110} (2003), 161--169.
%
\bibitem{MM3} Michalska, M.~--~ Mozgawa, W.: 
$\alpha$ -isoptics of a triangle and their connection to $\alpha$-isoptic of an oval, 
\textit{Rend. Semin. Mat. Univ. Padova}, Vol \textbf{133} (2015), p. 159--172
%
\bibitem{het} Miernowski, A.~--~ Mozgawa, W.: 
On some geometric condition for convexity of isoptics,
\textit{ Rend. Semin. Mat., Torino} \textbf{55}, No.2 (1997), 93--98.
%
\bibitem{Mi}
Milnor,~J.,
Curvatures of left Invariant metrics on Lie groups. 
{\it Advances in Math.,}  {\bf 21},  293--329 (1976).
%
\bibitem{M97}
{Moln{\'a}r,~E.,}
The projective interpretation of the eight 3-di\-men\-sional homogeneous geometries. 
\emph{Beitr. Algebra Geom.,}
{\bf38} No.~2, 261--288, (1997).
%
\bibitem{MoSzi10}
{Moln{\'a}r,~E.~--~Szil\'agyi,~B.,}
Translation curves and their spheres in homogeneous geometries.
\textit{Publ. Math. Debrecen,}
{\bf 78/2}, 327-346 (2010).
%
\bibitem{MSz}
{Moln{\'a}r,~E.~--~Szirmai,~J.,}
Symmetries in the 8 homogeneous 3-geometries.
\textit{Symmetry Cult. Sci.,}
{\bf 21/1-3}, 87-117 (2010).
%
\bibitem{MSz12}
{Moln{\'a}r,~E.~--~Szirmai,~J.,}
Classification of $\SOL$ lattices.
\textit{Geom. Dedicata,}
{\bf 161/1}, 251-275 (2012).
%
\bibitem{MSzV}
Moln{\'a}r,~E.~--~Szirmai,~J.~--~Vesnin,~A.,
Projective metric realizations of cone-manifolds with singularities along 2-bridge knots and links.
{\it J. Geom.,}  {\bf 95}, 91-133 (2009).
%
\bibitem{MSzV14}
Moln{\'a}r,~E.~--~Szirmai,~J.~--~Vesnin,~A.,
Packings by translation balls in $\SLR$. 
{\it J. Geom.,} {\bf 105}(2), 287--306 (2014)
%
\bibitem{tizenegy} Mozgawa, W.~--~ Skrzypiec, M.: 
Crofton formulas and convexity condition for secantopics, 
\textit{Bull. Belg. Math. Soc. - Simon Stevin} {\bf 16}, No. 3 (2009), 435--445.
%
\bibitem{Odehnal} Odehnal, B.: Equioptic curves of conic sections, 
\textit{J. Geom Graph.} \textbf{14} No.1 (2010), 29--43.
%
\bibitem{S}
Scott,~P.:
The geometries of 3-manifolds,
{\it Bull. London Math. Soc.}  {\bf 15} (1983), 401--487.
%
\bibitem{Si} Siebeck,~F.~H.~: 
\"Uber eine Gattung von Curven vierten Grades, welche mit den elliptischen Funktionen zusammenh\"angen,
\textit{J. Reine Angew. Math.} \textbf{57} (1860), 359--370; \textbf{59} (1861), 173--184.
%
\bibitem{tiz} Skrzypiec, M.: 
A note on secantopics, 
\textit{Beitr. Algebra Geom.} \textbf{49} No. 1 (2008), 205--215.
%
\bibitem{Sz13-1}
Szirmai,~J.,
A candidate to the densest packing with equal balls in the Thurston geometries. 
{\it Beitr. Algebra Geom.,} {\bf 55}(2), 441--452 (2014).
%
\bibitem{Sz18}
Szirmai,~J.,
Bisector surfaces and circumscribed spheres of tetrahedra derived
by translation curves in $\SOL$ geometry. 
{\it New York J. Math.,} {\bf 25}, 107--122 (2019).
%
\bibitem{Sz13-2}
Szirmai,~J.,
The densest translation ball packing by fundamental lattices in $\SOL$ space.  
{\it Beitr. Algebra Geom.,}  {\bf 51}(2)  353--373 (2010). 
%
\bibitem{Sz16}
Szirmai,~J.,
$\NIL$ geodesic triangles and their interior angle sums.  
{\it Bull. Braz. Math. Soc. (N.S.),} {\bf 49}  761--773 (2018), DOI: 10.1007/s00574-018-0077-9.
%
\bibitem{Sz20}
{Szirmai,~J.:}
Triangle angle sums related to translation curves in $\SOL$ geometry,
\textit{Stud. Univ. Babes-Bolyai Math.} {\bf 67} (2022), 621--631, arXiv: 1703.06646, doi: 10.24193/subbmath.2022.3.14.
%
\bibitem{Sz12}
{Szirmai,~J.:}
Lattice-like translation ball packings in $\NIL$ space.
\textit{Publ. Math. Debrecen}, {\bf 80/3-4} (2012), 427--440, DOI: 10.5486/PMD.2012.5117.
%
\bibitem{Sz202} 
Szirmai,~J.,
Interior angle sums of geodesic triangles in $\SXR$ and $\HXR$ geometries, 
{\it Bull. Academ. De Stiinte A Rep. Mol.}, {\bf{93}} Num 2 (2020), 44--61.
%
\bibitem{Sz22}
{Szirmai,~J.:}
Apollonius surfaces, circumscribed spheres of tetrahedra, Menelaus' and Ceva's theorems in $\SXR$ and $\HXR$ geometries,
\emph{Quarterly Journal of Mathematics}, {\bf 73} (2022), 477--494, doi: 10.1093/qmath/haab038, arXiv: 2012.06155.
%
\bibitem{Sz22-3}
{Szirmai,~J.:}
Classical Notions and Problems in Thurston Geometries,
\emph{International Electronic Journal of Geometry},
{\bf 16} No.2 (2023), 608--643, doi: 10.36890/IEJG.1221802, arXiv: 2203.05209.
%
\bibitem{Sz23}
{Szirmai,~J.:}
Fibre-like cylinders, their packings and coverings in $\SLR$ space,
\emph{Results Math.}, (2024), doi: 10.1007/s00025-024-02152-0, arXiv: 2306.05721.
%
\bibitem{T}
Thurston,~W.~P. (and Levy,~S. editor),
{\it Three-Dimensional Geometry and Topology}.  Princeton University Press,  Princeton, New Jersey, vol. {\bf 1} (1997).
%
\bibitem{yates} Yates,~R.~C.: A handbook on curves and their properties, 
J.~W.~Edwards, Ann. Arbor, (1947), 138--140. 
%
\bibitem{Wi} Wieleitener,~H.~: Spezielle ebene Kurven. Sammlung Schubert LVI,
\textit{G\"oschen'sche Verlagshandlung}. Leipzig, (1908).
%
\bibitem{Wu71-1} Wunderlich,~W.~: Kurven mit isoptischem Kreis,
\textit{Aequat. math}. \textbf{6} (1971), 71-81.
%
\bibitem{Wu71-2} Wunderlich,~W.~: Kurven mit isoptischer Ellipse,
\textit{Monatsh. Math.} \textbf{75} (1971), 346-362.
\end{thebibliography}
\end{document}